\documentclass[nosubfloats, final]{l4dc2025}

\usepackage{microtype}
\usepackage{tikz}     
\usepackage{booktabs} 
\usepackage{caption} 
\usepackage{enumitem} 
\usepackage{algorithm}
\usepackage{algorithmic}
\usepackage{multirow}
\usepackage{mathtools}
\usepackage{wrapfig}
\usepackage{hyperref}
\usetikzlibrary{calc}

\def\N{\mathbb{N}}
\def\R{\mathbb{R}}

\def\S{\mathbb{S}}

\def\col{\textrm{col}}
\def\Tini{T_{\textrm{ini}}}
\def\Up{U_{\textrm{p}}}
\def\Yp{Y_{\textrm{p}}}
\def\Wp{W_{\textrm{p}}}
\def\Uf{U_{\textrm{f}}}
\def\Yf{Y_{\textrm{f}}}
\def\Wf{W_{\textrm{f}}}
\def\Tf{T_{\textrm{f}}}
\def\ini{\textrm{ini}}
\def\uini{u^{\textrm{ini}}}
\def\yini{y^{\textrm{ini}}}
\def\wini{w^{\textrm{ini}}}

\def\lp{\left(}
\def\rp{\right)}

\def\ytil{\Tilde{y}}
\def\r{\textrm{ref}}

\def\wr{w^{\operatorname{ref}}}
\def\ftil{\Tilde{f}}
\newcommand{\tup}[1]{\textrm{#1}}
\DeclareMathOperator*{\argmin}{argmin} 

\title[DeePC-Hunt]{DeePC-Hunt: Data-enabled Predictive Control Hyperparameter Tuning via Differentiable Optimization}
\usepackage{times}

\coltauthor{
  \Name{Michael Cummins}{\normalfont\textsuperscript{1}} \Email{m.cummins24@imperial.ac.uk}\\
  \Name{Alberto Padoan}{\normalfont\textsuperscript{2}} \Email{apadoan@ece.ubc.ca}\\
  \Name{Keith Moffat}{\normalfont\textsuperscript{3}} \Email{kmoffat@ethz.ch}\\
  \Name{Florian D$\ddot{\text{o}}$rfler}{\normalfont\textsuperscript{3}}\Email{dorfler@control.ee.ethz.ch}\\
  \Name{John Lygeros}{\normalfont\textsuperscript{3}} \Email{jlygeros@ethz.ch}\\
  \addr {\normalfont\textsuperscript{1}}Department of Electrical and Electronic Engineering, Imperial College London\\{\normalfont\textsuperscript{2}}Department of Electrical and Computer Engineering, University of British Columbia\\{\normalfont\textsuperscript{3}}Department of Electrical Engineering and Information Technology, ETH Z$\ddot{\text{u}}$rich}
\begin{document}

\maketitle

\begin{abstract}
 This paper introduces \textbf{D}ata-\textbf{e}nabl\textbf{e}d \textbf{P}redictive \textbf{C}ontrol \textbf{H}yperparameter T\textbf{un}ing via Differen\textbf{t}iable Optimization (DeePC-Hunt), a backpropagation-based method for automatic hyperparameter tuning of the DeePC algorithm. The necessity for such a method arises from the importance of hyperparameter selection to achieve satisfactory closed-loop DeePC performance. The standard methods for hyperparameter selection are to either optimize the open-loop performance, or use manual guess-and-check. Optimizing the open-loop performance can result in unacceptable closed-loop behavior, while manual guess-and-check can pose safety challenges. DeePC-Hunt provides an alternative method for hyperparameter tuning which uses an approximate model of the system dynamics and backpropagation to directly optimize hyperparameters for the closed-loop DeePC performance. Numerical simulations demonstrate the effectiveness of DeePC in combination with DeePC-Hunt in a complex stabilization task for a nonlinear system and its superiority over model-based control strategies in terms of robustness to model misspecifications.\footnote{GitHub repo for reproducing experiments: 
 \href{https://github.com/michael-cummins/DeePC-HUNT}{https://github.com/michael-cummins/DeePC-HUNT}}
\end{abstract}


\section{Introduction}
Over the past decade, direct data-driven control methods have experienced a surge of interest~\citep{coulson2019dataenabled,persis, vanwaarde2020data_informativity,xue2021datadriven_SLS},  primarily fueled by the increasing adoption of machine learning techniques and the availability of vast datasets. Unlike traditional control design paradigms that rely on system identification followed by model-based
control, these methods compute control actions directly from data. Leveraging tools from behavioral system theory, 
uncertainty quantification 
and convex optimization, data-driven methods have shown promise in simplifying the control design process, while achieving satisfactory control performance with reduced implementation effort.

Among the methods in this expanding literature, Data-enabled Predictive Control (DeePC) \citep{coulson2019dataenabled} stands out as an effective direct data-driven control algorithm.
Operating similarly to Model Predictive Control (MPC) in a receding-horizon manner \citep{borelliBook}, DeePC circumvents the need for an accurate state-space model, relying instead on a Hankel matrix derived from offline input/output raw data. Although stability guarantees have been established only under specific assumptions \citep{coulson2019dataenabled}, the DeePC algorithm has demonstrated remarkable performance in controlling nonlinear systems affected by noise. This performance is mainly attributed to regularization techniques \citep{coulson2019regularized} and is theoretically justified by tools from distributionally robust optimization \citep{distributionally_robust_deepc}. However, the performance of DeePC is often  sensitive to the regularization parameters \citep{dörfler2021bridging}, presenting challenges when experiments are difficult, costly, or unsafe. 

Today, the two standard methods for regularization tuning are analytical methods for the open-loop and manual guess-and-check. Open-loop methods, such as Final Control Error \citep{chiuso2023harnessing}, rely on structural assumptions of the objective function and often yield overly conservative performance, as they do not account for DeePC's ability to replan at each timestep.  Manual guess-and-check methods, on the other hand, optimize the closed-loop performance but require experiments on the real system, which may be costly or even impossible. In contrast, DeePC-Hunt directly optimizes hyperparameters for receding-horizon closed-loop performance via offline backpropagation, without requiring strong modeling assumptions or direct interaction with the system.

\begin{wrapfigure}{r}{0.43\linewidth}
    \vspace{-10pt}
        \includegraphics[width=\linewidth]{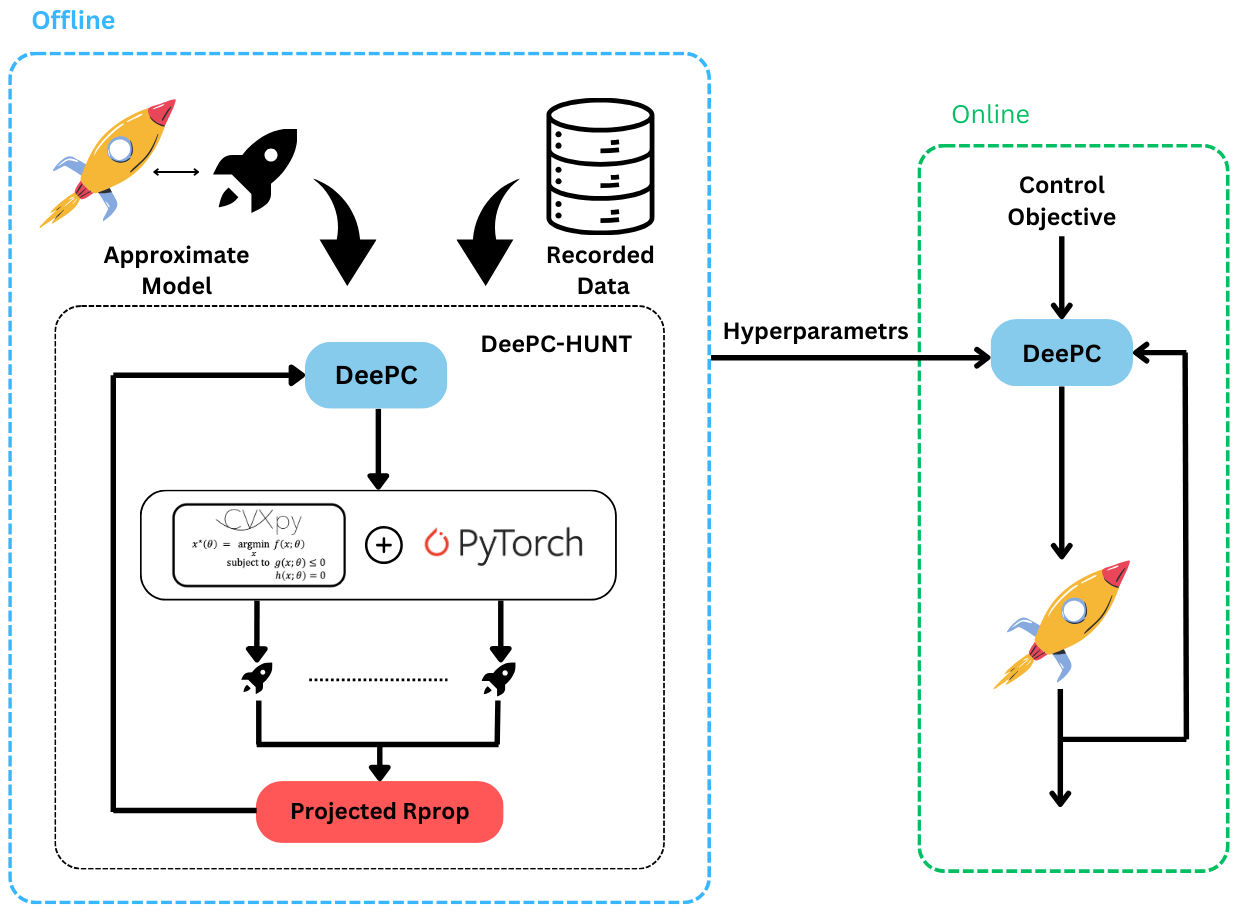}
        \caption{
        The DeePC policy is instantiated using CvxpyLayers and PyTorch to enable automatic differentiation. Simulations are carried out on an approximate model of the system and projected resilient backpropagation is used to update the hyperparameters.}
        \captionsetup{skip=0pt}
    \vspace{-0.1cm}
\end{wrapfigure}

This paper introduces \textbf{D}ata-\textbf{e}nabl\textbf{e}d \textbf{P}redictive \textbf{C}ontrol \textbf{H}yperparameter T\textbf{un}ing via Differen\textbf{t}iable Optimization (DeePC-Hunt), a backpropagation-based method designed to automate hyperparameter tuning for DeePC. \emph{DeePC-Hunt interprets DeePC, i.e., the solution of the DeePC optimization, as a control policy, and interprets the DeePC regularization hyperparameters as the parameters of the policy}.
Leveraging this DeePC-as-policy interpretation and an approximate model of the system dynamics, DeePC-Hunt optimizes the hyperparameters (to local optimality) by backpropagating \citep{backprop} closed-loop performance.
Specifically, DeePC-Hunt implements a (constrained) variant of the resilient backpropagation algorithm \citep{Rprop} and widely-used automatic differentiation tools~\citep{Torch} to directly optimize the regularization hyperparameters of a DeePC policy based on the closed-loop performance of the policy deployed on the approximate model in offline experiments. 
Our findings suggest that, given a reasonably approximated model, DeePC-Hunt yields performant hyperparameters for the true system. Moreover, DeePC-Hunt demonstrates superior robustness to model inaccuracies over traditional model-based control methods. 

\textbf{Related work:} DeePC-Hunt is motivated and inspired by policy optimization algorithms \citep{policy_optimisation}, which iteratively refine control policy parameters to minimize cumulative cost using different versions of the (projected) gradient descent algorithm. Unlike MPC or DeePC, DeePC-Hunt yields policies expressed as explicitly differentiable functions mapping states to actions. The implicit function theorem \citep{amos2021optnet} can be then used to differentiate the solution map of a Convex Optimization Control Policy (COCP), where Karush-Kuhn-Tucker (KKT) conditions are differentiated to compute the gradient of the control input with respect to the hyperparameters. The work closest to ours are \citep{amos2019differentiable, zuliani2024bpmpc}, which extend the theory of differentiating KKT conditions theory to MPC policies for automatic parameter tuning via gradient descent methods. However, \citep{amos2019differentiable, zuliani2024bpmpc} explicitly \ rely on state-space models, which can lead to unacceptable performance if the model is inaccurate. A broader approach to policy optimization for COCPs is presented by \citep{agrawal2019learning}, with experimental validation demonstrating the effectiveness of differentiable optimization for tuning a wide range of COCPs. The use of an approximate model in DeePC-Hunt is further motivated by results in sample-efficient policy optimization \citep{song2024learning, wiedemann2023training}, where policy gradients are computed directly rather than estimated. Unlike previous approaches, DeePC-Hunt optimizes DeePC policies by integrating data from a system with data obtained from numerical simulations of an approximate model of the system, resulting in improved closed-loop performance.

\textbf{Contributions:} The main contributions of the paper are threefold. 
\begin{enumerate}[label=\roman*), ref=\roman*, itemsep=0pt, topsep=0pt, parsep=0pt, partopsep=0pt]
\item We formulate the problem of hyperparameter tuning for DeePC policies, ensuring compatibility with backpropagation, and propose the DeePC-Hunt algorithm as an effective solution. 
\item To efficiently implement DeePC-Hunt, we introduce a variant of the resilient backpropagation algorithm \citep{Rprop} that supports the use of constraints (e.g. box constraints) using off-the-shelf automatic differentiation tools~\citep{Torch}. 
\item We validate DeePC-Hunt on a challenging benchmark task: landing a Vertical Takeoff Vertical Landing (VTVL) vehicle on an oceanic platform using a 
 high-fidelity Gym simulation environment~\citep{brockman2016openai}.  
\end{enumerate}

\textbf{Paper organization:} 
The paper is organized as follows. Section \ref{sec:background} presents a concise overview of differentiable convex optimization layers, DeePC, and a version of Resilient Backpropagation~\citep{Rprop}. Section \ref{sec:deepcHunt} introduces the  DeePC-Hunt algorithm. Section \ref{sec:numerical} demonstrates the robustness and performance of DeePC-Hunt on the VTVL landing task, comparing its performance against MPC under model mismatch. Section \ref{sec:conclusion} concludes the paper with a summary of our findings.

\textbf{Notation:}
The set of real numbers is denoted by $\R$. 
The set of real $n$-dimensional vectors is denoted by $\R^{n}$.
The set of real $n \times m$-dimensional matrices is denoted by $\R^{n \times m}$.
The set of positive integers is denoted by $\N$. The set of non-negative real numbers is denoted by $\R_+$. The set of  ${n\times n}$-dimensional symmetric positive definite matrices  is denoted by ${\mathbb{S}_{+}^{n \times n}}$. The set of closed, convex cones in $\R^m$ is denoted by $\mathcal{C}^m$. The transpose of the matrix ${M \in \R^{p \times m}}$ is denoted by 
$M^{\intercal}$. The $p$-norm of the vector ${x\in \R^n}$ is denoted by $|x|_p$. The weighted norm of the vector ${x\in \R^n}$ induced by the matrix ${Q \in \S_+^{n \times n}}$ is denoted by ${|x|_Q}$  and defined as ${|x|_Q= (x^\intercal Q x)^{1/2}}$.\
The $\text{i}^{\text{th}}$ entry of a vector ${v^{\operatorname{d}} \in \R^n}$ is denoted by $v^{\operatorname{d}|i}$. The expectation operator is denoted by $\mathbb{E}$. The function ${\operatorname{sign}:\R \to \{-1,0,1\}}$ is defined as $-1$ if the argument is negative, $1$ if the argument is positive, and $0$ if the argument is zero,  respectively. Given two vectors $v\in \R^{n}$ and $w\in \R^{m}$, we define $\col(v,w) := \left(v^\intercal,w^\intercal\right)^\intercal$. Given $T\in\N$, the {Hankel matrix} of depth ${L\in\N}$, with ${L \le T}$,  associated with the vector ${w \in \R^{qT}}$ 
is denoted by $H_L(w)$ (see, e.g.,~\citep{coulson2019dataenabled} for the definition).

\section{Background}\label{sec:background}

\subsection{Differentiable Convex Optimization Layers}\label{sec:DCOL}
Consider the parameterized convex optimization problem
 \begin{align}\label{eq:general_opt}
     \min_{x \in \R^n}\:\: f(x | \theta) \:\: \text{s.t.} \:\: g(x|\theta) \leq 0, \:\: h(x|\theta) = 0,
 \end{align}
\noindent
where $f: \R^n \to \R$ and $g: \R^n \to \R^{m_1}$ are convex functions, $h: \R^n \to \R^{m_2}$ is affine and  ${\theta \in \R^{\ell}}$ represents the vector of parameters defininng $f$, $g$ and $h$, respectively.
The parameter $\theta$ is assumed to belong to a given set ${\Theta\subseteq \R^\ell}$, which represents the problem parameters.  The solution to the class of convex optimization problems of the form~\eqref{eq:general_opt}, known as \emph{disciplined parameterized programs} \citep{layers}, may be defined as a mapping ${\theta \mapsto x^*(\theta)}$ from the parameter $\theta$ to the (global) solution $x^*(\theta)$ of~\eqref{eq:general_opt}. Formally, we define the solution map of \eqref{eq:general_opt} as
\begin{equation}\label{eq:general_sol}
    s(\theta) := \theta \mapsto x^*(\theta). 
\end{equation} 
\noindent    
Differentiable convex optimization layers \citep{layers} provides a method for differentiating the  solution-map \eqref{eq:general_sol} with respect to the parameter $\theta$. To this end, \eqref{eq:general_opt} is first transformed into a convex cone problem defined as 
\begin{align}
\label{eq:cone_program}
    \min_{(x,\nu) \in \R^n \times \R^m} \:\: c^\intercal x 
    \:\: \text{s.t.} \:\: Ax + \nu = b, 
    \:\: (x,\nu) \in \mathbb{R}^n \times \mathcal{K}, 
\end{align}  

\noindent
where ${A \in \R^{m\times n}}$, ${b \in \R^m}$, ${c \in \R^n}$ and ${\mathcal{K} \in \mathcal{C}^m}$ are transformations of the problem parameter~$\theta$. Such a transformation is typically performed through domain-specific languages for convex programming, such as CVXPY 
\citep{diamond2016cvxpy}, where the subset of parameters $(A,b,c)$ that affect the gradient are determined by a sparse (linear) projection $L(\theta) = (A,b,c)$ \citep{rewriting}. 
Thus, the solution map \eqref{eq:general_sol} is decomposed as
\begin{equation}\label{eq:cone_solution_map}
    s(\theta) = (R \circ S \circ L)(\theta),
\end{equation}
where ${L: \R^\ell \to \R^{m \times n} \times \mathbb{R}^m \times \mathbb{R}^n}$ maps the problem data of \eqref{eq:general_opt} to the problem data of \eqref{eq:cone_program},
${S : \R^{m \times n} \times \mathbb{R}^m \times \mathbb{R}^n \to \R^n \times \mathcal{K}}$ is the solution map of \eqref{eq:cone_program} and ${R: \R^n \times \mathcal{K} \to \R^n}$ maps the solution of \eqref{eq:cone_program} to the solution of \eqref{eq:general_opt}. Specifically, ${S(A,b,c) = (x^*,\nu^*)}$ and $R(x^*,\nu^*) = x^*$. The gradient of the solution map is then ${\nabla s(\theta) = \nabla L(\theta)\nabla S(A,b,c)\nabla R(x^*,\nu^*)}$. Since $R$ and $L$ are linear mappings by design, calculating their gradient is elementary. 
Furthermore, the gradient $\nabla S(A,b,c)$ may be obtained by differentiating the KKT conditions of \eqref{eq:cone_program}, as detailed in \citep{agrawal2020differentiating} (see also \citep{rewriting, layers}). 

The CvxpyLayers \citep{layers} Python package integrates CVXPY and PyTorch \citep{Torch} to automate the transformation \eqref{eq:general_opt} $\mapsto$ \eqref{eq:cone_program} and the differentiation of \eqref{eq:cone_solution_map} with respect to $\theta$. DeePC, explained in the following section, can be interpreted as a convex optimization problem of the form \eqref{eq:general_opt}. CvxpyLayers will be later used in the implementation of DeePC-Hunt to instantiate the DeePC policy \eqref{eq:DeePC}, allowing us to differentiate the control inputs with respect to the regularization parameters.

\subsection{The DeePC Algorithm}\label{sec:DeePC}
Consider a discrete-time dynamical system described by the equations
\begin{align}\label{eq:model} 
    x_{k+1} = f\lp x_k,u_k\rp, \:\: y_k = h\lp x_k,u_k\rp, 
\end{align}
\noindent
where ${u\in\R^m}$, ${x\in \R^n}$ and ${y\in\R^p}$ represent the input, state, and output of system \eqref{eq:model}, respectively. 
The DeePC algorithm leverages raw data matrices derived from offline input/output measurements of system \eqref{eq:model} as a predictive model. Assuming the availability of input/output data recorded offline from system \eqref{eq:model}, let ${u^{\tup d} = \col(u^{\tup d | 1},\dots,u^{\tup d | T}) \in \R^{mT}}$ and ${y^{\tup d} = \col(y^{\tup d | 1},\dots,y^{\tup d | T}) \in \R^{pT}}$ be vectors containing an input sequence of length ${T\in \N}$ applied to system~\eqref{eq:model} and the corresponding output sequence, respectively. Let ${q=m+p}$ and define the data vector 
$w^{\text{d}} = \col(u^{\tup d}, y^{\tup d})\in\R^{qT}$. 
For given initial and future time horizons $\Tini \in \N$ and $\Tf \in \N$, we define the input and output data matrices for the past and future timesteps as
\scalebox{0.75}{$\begin{pmatrix}
    \Up \\ \Uf 
\end{pmatrix}$} = $H_{\Tini+\Tf}(u^{\textup{d}})$ and \scalebox{0.75}{$\begin{pmatrix}
    \Yp \\ \Yf
\end{pmatrix}$} = $H_{\Tini + \Tf}(y^{\textup{d}})$, respectively.
Furthermore, we define the past and future combined input-and-output data matrices as 
$\Wp = $\scalebox{0.75}{$\begin{pmatrix}
    \Up \\ \Yp 
\end{pmatrix}$} and $\Wf = $\scalebox{0.75}{$\begin{pmatrix}
    \Uf \\ \Yf 
\end{pmatrix}$}, respectively. 

Given a reference trajectory $\wr=\col(u^{\r|1},\dots,u^{\r|\Tf},y^{\r|1},\dots,y^{\r|\Tf})\in \R^{q\Tf}$, an input constraint set $\mathcal{U}\subseteq \R^{m\Tf}$, an output constraint set $\mathcal{Y} \subseteq \R^{p\Tf}$, an input cost matrix $R \in \mathbb{S}_+^{m\times m}$, an output cost matrix $Q\in \mathbb{S}_+^{p\times p}$, a regularization function $\psi:\R^{T-\Tini-\Tf+1} \times \R^{p\Tini} \to \R^r_+$, and a weight vector $\lambda \in \R^r_+$, the DeePC algorithm operates in a receding-horizon fashion by iteratively solving the optimization problem
\begin{align}
\label{eq:DeePC}
    \underset{u,y,g,\sigma_y}{\text{min}}\,
    &  
    \displaystyle
    \sum_{i=1}^{\Tf} |y_i - y^{\r|i}|_Q^2 + |u_i - u^{\r|i}|_R^2
    +\lambda^\intercal\psi(g,\sigma_y),  \nonumber \\
    \text{s.t.\,}
    & \begin{pmatrix}
    \Up \\ \Yp \\ \Uf \\ \Yf 
    \end{pmatrix}g
    =\begin{pmatrix}
    \uini_k \\ \yini_k \\ u \\ y
    \end{pmatrix} 
    + \begin{pmatrix}
        0\\
        \sigma_y\\
        0\\
        0
    \end{pmatrix}, \:\: 
    (u,y)\in\mathcal{U}\times \mathcal{Y},
\end{align}   
where \(\uini_k = \operatorname{col}(u_{k-\Tini}, \dots, u_{k-1}) \in \mathbb{R}^{m\Tini}\) and \(\yini_k = \operatorname{col}(y_{k-\Tini}, \dots, y_{k-1}) \in \mathbb{R}^{p\Tini}\) are vectors containing the \(\Tini\) most recent input and output measurements at time \(k\), respectively.

Let \(\wini_k = \operatorname{col}(\uini_k, \yini_k) \in \mathbb{R}^{q\Tini}\) and note that \eqref{eq:DeePC} is a parametric convex optimization problem with parameter $\lambda$, we define its solution map in terms of $\lambda$, ${\wini_k}$ and ${\wr_k}$ as $\mathcal{S}\left(\lambda,\wini_k,\wr\right) := \left(\lambda,\wini_k,\wr\right) \mapsto (w^*,g^*,\sigma_y^*),$ where ${w^* = \col(u^*,y^*)}$ and ${(w^*,g^*,\sigma_y^*)}$ is the optimal solution of the optimization problem~\eqref{eq:DeePC}. Furthermore, defining the (linear) projection $\mathcal{P}(w^*,g^*,\sigma_y^*) = u^*_0,$ where $u_0^* \in \R^m$ is the vector defined by the first $m$ entries of $w^*$, we rewrite the DeePC control policy as 
\begin{equation}\label{eq:deepc_policy}
    \pi^\lambda(\wini_k,\wr) := (\mathcal{P} \circ \mathcal{S})\left(\lambda,\wini_k,\wr\right).
\end{equation}
The original formulation of the DeePC algorithm in \citep{coulson2019dataenabled} employs a one-norm regularizer on the decision variable ${g \in \R^{T-\Tini-\Tf+1}}$ to promote the selection of low-complexity models. 
Alternatively, as detailed in~\citep{dörfler2021bridging}, one may use the Elastic-Net Regularization function which leads to consistent predictions and connects to classic Subspace Predictive Control (SPC) algorithms ~\citep{SPC}. In this work, we combine such regularizer with a one-norm regularization to promote sparsity in the slack variable $\sigma_y$ using the regularization function
\begin{equation}\label{eq:regulariser}
    \lambda^\intercal \psi(g,\sigma_y) = \lambda^0 \left| (I-\Pi)g \right|_2^2 + \lambda^1 \left| g \right|_1 + \lambda^2|\sigma_y|_1,
\end{equation} 
where ${\lambda = \col(\lambda^0,\lambda^1,\lambda^2)}$ 
and $\Pi = $
\scalebox{0.75}{$
\begin{pmatrix}
\Wp \\ \Yf 
\end{pmatrix}^{\dagger}
\begin{pmatrix}
\Wp \\ \Yf 
\end{pmatrix}
$}.

\subsection{Projected Resilient Backpropagation}\label{sec:Rprop}
Consider the constrained optimization problem
    $$\min_{\lambda\in\Lambda} \:\phi(\lambda),$$
where ${\Lambda \subseteq \R^r}$ and ${\phi: \R^r \to \R}$ is differentiable.
Finding a locally optimal solution ${\lambda^{*} \in \R^r}$ may be challenging 
with a first-order iterative optimization algorithm 
if $|\nabla \phi(\lambda)|$ is small for all $\lambda \in \R^n$, regardless of the proximity to a stationary point. 
Resilient backpropagation \citep{Rprop} is a heuristic first-order iterative optimization algorithm that mitigates this issue by optimizing over each scalar element of the decision variable ${\lambda = \col(\lambda^0,\dots,\lambda^{r-1})}$ using only the \textit{sign} of the gradient and an adaptive step-size. 
The update rule of the algorithm is defined as
\noindent
\begin{align}\label{eq:rprop}
    \lambda_{k+1}^i = \lambda_{k}^i - \eta_{k}^i \textrm{sign} (\nabla_{\lambda^i}\phi(\lambda_{k}^i)),\:\:
    \eta_{k+1}^i = \text{Rprop}(\eta_{k}^i,\lambda_{k}^i,\lambda_{k+1}^i|\eta^{\text{max}}, \eta^{\text{min}}, \beta, \alpha), 
\end{align}
where ${\operatorname{Rprop}:\R_+ \times \R \times \R \to \R}$ is the step-size generator defined as
\begin{align*}
    \text{Rprop}(\eta_{k}^i,\lambda_{k}^i,\lambda_{k+1}^i|\eta^{\text{max}}, \eta^{\text{min}}, \beta, \alpha) &= \begin{cases}
        \min\{\alpha \eta_{k}^i , \eta_{\max} \}, & \textrm{if} \: \nabla_{\lambda^i}\phi(\lambda_{k+1}^i) \nabla_{\lambda^i}\phi(\lambda_{k}^i) < 0,\nonumber  \\ 
        \max\{\beta\eta_{k}^i , \eta_{\min} \}, & \textrm{if} \: \nabla_{\lambda^i}\phi(\lambda_{k+1}^i) \nabla_{\lambda^i}\phi(\lambda_{k}^i) > 0,\nonumber \\
        \eta_{k}^i, & \textrm{if} \: \nabla_{\lambda^i}\phi(\lambda_{k+1}^i) \nabla_{\lambda^i}\phi(\lambda_{k}^i) = 0, \nonumber
    \end{cases}
\end{align*}

\noindent
with ${(\eta^{\text{max}}, \eta^{\text{min}}, \beta, \alpha) \in \R_+^4}$ being the design parameters representing the minimum step-size, maximum step-size, decay factor, and growth factor, respectively. 

Resilient backpropagation is particularly useful in the context of DeePC-Hunt, because the (locally) optimal DeePC regularization parameters can be orders of magnitude away from an initial guess~\citep{dörfler2021bridging}. 
However, each regularization term $\lambda^i$ must be non-negative. 
To enforce this constraint, we propose a variant of \eqref{eq:rprop} defined by the update rule
\begin{align}\label{eq:projected_rprop}
    \gamma_{k+1}^i =& 
    \gamma_{k}^i - \eta_{k}^i \textrm{sign}(\nabla_{\lambda^i}\phi(\lambda_{k}^i)), \:\:
     \eta_{k+1}^i = \text{Rprop}(\eta_{k}^i,\gamma_{k}^i,\gamma_{k+1}^i|\eta^{\text{max}}, \eta^{\text{min}}, \beta, \alpha), \quad \forall i = 1,\cdots,r\nonumber \\
     \lambda_{k+1} =& P_\Lambda \left( \gamma_{k+1} \right),
\end{align}
where ${P_\Lambda}$ is the projection operator onto the set $\Lambda$,  defined as ${P_\Lambda(\lambda) = \argmin_{{\mu\in \Lambda}} \frac{1}{2}|\lambda-\mu|^2}$.  The projection onto the set $\Lambda$ ensures that the constraint $\lambda^k \in \Lambda$ is satisfied at every iteration. To the best of the authors' knowledge, this is the first work to incorporate resilient backpropagation with a projection scheme to enforce feasibility in the iterates of \eqref{eq:rprop}. Throughout this paper, we assume that $P_\Lambda (\lambda)$ is easy to compute, as is the case for a range of common sets~\citep{parikh2014proximal}, including half-spaces, hyperrectangles, simplices, and certain cones.

\section{DeePC-Hunt}\label{sec:deepcHunt}

The regularization parameter $\lambda$ has a significant impact on the closed-loop performance of DeePC \citep{coulson2019dataenabled}. DeePC-Hunt tunes $\lambda$ to optimize the closed-loop performance of a DeePC policy implemented on a simplified or inaccurate surrogate model,  leveraging the intuition that optimal regularization parameters for ``similar'' systems are typically ``close.'' This approach enables offline hyperparameter tuning with a minimal sim-to-real gap:  by optimizing $\lambda$ with respect to the closed-loop performance on the surrogate model, we obtain a regularizer $\lambda$ that leads to acceptable performance on the true system.  
We consider a surrogate model described by the equations
\begin{align}\label{eq:model_approximate} 
    \tilde{x}_{k+1} = \tilde{f}\lp \tilde{x}_k,\tilde{u}_k\rp, \:\: \tilde{y}_k = \tilde{h}\lp \tilde{x}_k,\tilde{u}_k\rp, 
\end{align}
where $\tilde{u}\in\R^m$, $\tilde{x}\in \R^n$ and $\tilde{y}\in\R^p$, respectively. The functions $\tilde{f}$ and $\tilde{h}$ are assumed to approximate the input-output behavior of system~\eqref{eq:model} with respect to a given metric (e.g.,  $\mathcal{L}_2$ gain). 
The approximate model may be derived from the physics of the system~\eqref{eq:model} or through the use of model reduction techniques (see, e.g.,~\citep{antoulas2005approximation} for an overview of available techniques).  

Given the approximate model~\eqref{eq:model_approximate}, the vector $\wini_1$ consisting of the most recent $\Tini$ measurements at timestep 1, the reference trajectory $\wr$, and the DeePC policy $\pi^{\lambda}$ (Section \ref{sec:DeePC}), we define the approximate closed-loop cost, ${\tilde{C}^{\pi^{\lambda}}:\R^{q\Tini} \times \R^{q\Tf} \to \R_+}$, over the simulation length $N$ as
\begin{equation}\label{eq:closed_loop_cost}
    \tilde{C}^{\pi^{\lambda}}\left(\wini_1,\wr\right) = \sum_{i=1}^{N} |\ytil_i-y^{\r|i}|_Q^2 
    + |u^*_i-u^{\r|i}|_R^2,
\end{equation} 
where $u^*_i$ is the input produced by $\pi^\lambda\left(\wini_i,\wr\right)$ in \eqref{eq:deepc_policy}
and $\ytil_i$ is the output of the surrogate model given by \eqref{eq:model_approximate} when $u^*_i$ is applied.
At each timestep $i$, $\wini_{i+1}$ is given by inserting $u^*_i$ and $\ytil_i$ into $\wini_i$, and removing the measurements furthest in the past corresponding to time $i-\Tini$.
The initial state $\tilde{x}_0$ required by the surrogate model \eqref{eq:model_approximate} is determined from $\wini_1$ using a state estimation method (e.g., a Kalman Filter). We define $\zeta(\cdot)$ as the selection operator which implements the Kalman Filter, i.e., $\zeta(\wini_i) = \tilde{x}_{i-1}$. Moreover, we define the operators $\Omega_u(\cdot)$ and $\Omega_y(\cdot)$ such that for a given $\wini_i$, we have $\Omega_u(\wini_i) = (\uini_{i-\Tini+2},\dots,\uini_{i})$ and $\Omega_y(\wini_i) = (\yini_{i-\Tini+2},\dots,\yini_{i})$.

The aim of DeePC-Hunt is to find a $\lambda$ that is effective for all initial conditions. Therefore, the expected closed-loop cost over a probability distribution of initial conditions $\wini_1$ is optimized~\footnote{A similar procedure could be followed for $\wr$, 
if the DeePC policy needs to work for a distribution of reference trajectories as well.}. For simplicity, the $\wini_1$ distribution is taken to be the uniform distribution over the columns of $\Wp$, denoted $\mathcal{D}(\Wp)$. Thus, the expected closed-loop cost is
$
    \mathbb{E}_{\wini_1 \sim \mathcal{D}(\Wp)}\left[\tilde{C}^{\pi^{\lambda}}\left(\wini_1,\wr\right)\right]
$
and the constrained, non-convex DeePC-Hunt optimization problem is defined by \eqref{eq:opt_framework} the bilevel optimization problem 
\begin{align}\label{eq:opt_framework}
    \underset{\lambda \in \Lambda}{\min} & \:\:\mathbb{E}_{\wini_1\sim \mathcal{D}(\Wp)}\left[\tilde{C}^{\pi^{\lambda}}\left(\wini_1,\wr\right)\right], \\
    \text{s.t.} \nonumber
    &\:\: \tilde{x}_0 = \zeta (\wini_1), 
    \:\:\tilde{x}_{i+1} = \tilde{f}\lp \tilde{x}_i,\tilde{u}_i\rp, \:\: 
    \tilde{y}_i = \tilde{h}\lp \tilde{x}_i,\tilde{u}_i\rp, \nonumber  \\
    &\:\: \tilde{u}_i = \pi^\lambda\left(\wini_i,\wr\right), \:\:\wini_{i+1} = \col(\Omega_u(\wini_i), \tilde{u}_{i+1}, \Omega_y(\wini_i), \tilde{y}_{i+1}), \nonumber 
\end{align} 
where $\Lambda \subseteq \R_+^r$ is selected to enforce non-negative box constraints on each $\lambda^i$. Recalling \eqref{eq:deepc_policy}, $\pi^\lambda$ is the solution map of a disciplined parameterized program with a linear projection applied. Therefore, $\pi^\lambda$ may be differentiated using CvxpyLayers, which leverages the technique discussed in Section \ref{sec:DCOL}.
Thus, the projected resilient backpropagation method described in \eqref{eq:projected_rprop} can be used to compute a locally optimal solution $\lambda^*$ of \eqref{eq:opt_framework}.

The DeePC-Hunt optimization problem \eqref{eq:opt_framework} is an empirical risk minimization with dataset $\mathcal{D}(\Wp)$. Thus, while the expectation \eqref{eq:opt_framework} can be computed, it may be prohibitively expensive if $\mathcal{D}(\Wp)$ has many datapoints as the exact expectation requires numerous simulations to produce a single gradient step. To compute the approximate gradient of the expectation in \eqref{eq:opt_framework} in a computationally tractable manner, a similar approach is taken to Stochastic Gradient Descent (SGD) methods by computing an unbiased estimate of \eqref{eq:opt_framework} in the form of Monte Carlo samples, which is common practice when performing empirical risk minimization with large datasets \citep{LSCO}. Pseudo-code for the sample-estimate/Monte Carlo DeePC-Hunt implementation for solving \eqref{eq:opt_framework} is given in \algorithmref{alg:deepc}, where $B \in \N$ is the number of Monte Carlo samples taken at each step $k$ and $N_{\text{iter}} \in \N$ is the number of steps.

\begin{algorithm}[H]
    \caption{DeePC-HUNT}
    \label{alg:deepc}
    \begin{algorithmic}
    \STATE \textbf{Initialize:} $\eta^i_0$, $\lambda^i_0$, $\gamma^i_0$, $\forall i = 1,\dots,r$
    \FOR{$k = 1:N_{\text{iter}}$}
        \FOR{$j = 1:B$ \textbf{in parallel}}
            \STATE $w^{\ini,j}_1 \sim \mathcal{D}(\Wp)$
            \STATE ${\tilde{x}_0^j = \zeta (w^{\ini,j}_1)}$
            \FOR{$i = 0:N-1$}
                \STATE $\tilde{u}_i^{j} \leftarrow \pi^{\lambda^k} \left(w^{\ini,j}_i,\wr\right)$
                \STATE $(\tilde{x}_{i+1}^j, \ytil_{i}^{j}) \leftarrow 
                (\ftil (\tilde{x}_i^j, \tilde{u}_i^{j}), \tilde{h}( \tilde{x}_{i}^j, \tilde{u}_i^{j}))$
                \STATE $w^{\ini,j}_{i+1} \leftarrow \col(\Omega_u(w^{\ini,j}_i), \tilde{u}_{i}^j, \Omega_y(w^{\ini,j}_i), \tilde{y}_{i}^j)$
            \ENDFOR
        \ENDFOR
        \STATE $J(\lambda^k) \leftarrow \frac{1}{b}\sum_{j=1}^b \tilde{C}^{\pi^{\lambda}}\left(w^{\ini,j},\wr\right)$
        \FOR{$i = 1:r$ \textbf{in parallel}}
            \STATE $\gamma_{k+1}^i \leftarrow \gamma_{k}^i - \eta_{k}^i \textrm{sign} (\nabla_{\lambda} J(\lambda_{k}^i))$
            \STATE $\eta_{k+1}^i \leftarrow \text{Rprop} \left(\eta_{k}^i, \gamma_{k}^i,\gamma_{k+1}^i\:|\eta^{\text{max}},\eta^{\text{min}},\beta,\alpha\right)$
        \ENDFOR 
        \STATE $\lambda_{k+1} \leftarrow P_\Lambda (\gamma_{k+1})$ 
    \ENDFOR
    \end{algorithmic}
\end{algorithm}

\section{Numerical Simulation: landing a VTVL vehicle}\label{sec:numerical}
\subsection{Rocket Lander Gym Environment}

The performance of DeePC-Hunt is demonstrated via the task of safely landing a Vertical Takeoff and Vertical Landing (VTVL) rocket. The objective is to land the VTVL vehicle on a designated oceanic platform in a vertical position. The problem is introduced in \citep{Ferrante2017ARC}  and implemented using OpenAI's Gym Library~\citep{brockman2016openai}.

\begin{wrapfigure}{r}{0.5\linewidth}
\vspace{-25pt}
    \includegraphics[width=\linewidth]{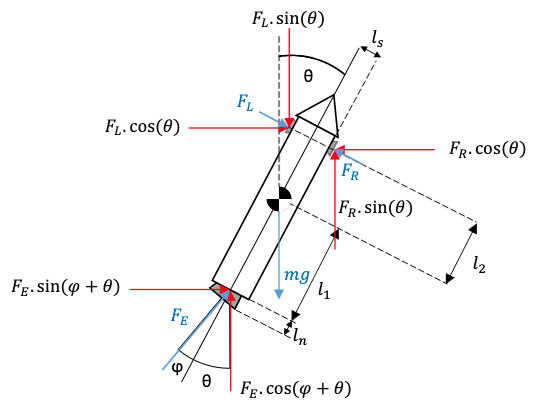}
  \caption{Free body diagram of the VTVL rocket from \citep{Ferrante2017ARC}.}
  \label{fig:free_body}
  \vspace{-30pt}
\end{wrapfigure}

The dynamics of the VTVL vehicle are modeled by the nonlinear equations of motion 
\begin{align}\label{eq:dynamics} 
    m\ddot{x} &= F_s\cos{(\theta)} - F_E\sin{(\varphi + \theta)}l_1 ,\\
    m\ddot{y} &= F_s\cos{(\theta)} - F_E\sin{(\varphi + \theta)}l_1 - mg , \nonumber\\
    m\ddot{\theta} &= F_E\sin(2\pi - \varphi)l_1 - F_sl_2\nonumber ,
\end{align}

\noindent
where $x(t)\in\mathbb{R}$ is the horizontal position of the rocket,
$y(t)\in\mathbb{R}$ is the vertical position of the rocket,
$\theta(t)\in\mathbb{R}$ is the vertical pitch of the rocket, 
$F_s(t)\in \mathbb{R}$ is the force exerted by the side engines, 
$F_E(t) \in \R_+$ is the force exerted by the main engine, 
$\varphi(t)\in\mathbb{R}$ is the heading angle of the main engine, 
${m(t)\in\mathbb{R}}$ is the time-varying mass of the rocket due to loss of fuel, ${l_1\in\mathbb{R}_+}$ is the length of the portion of the rocket below the center of gravity, and ${l_2\in\mathbb{R}_+}$ is the length of the portion of the rocket above the center of gravity, respectively. 
For simplicity, we assume that the landing pad is stationary, perfect state information is available, and control inputs are chosen as ${u = (F_E, F_s, \varphi) \in \mathbb{R}^3}$, thus facilitating full actuation. The state and input are also required to satisfy given state and input box constraints ${x(t)\in X}$, ${y(t)\in Y}$, ${\theta(t)\in \Theta}$, and ${u(t)\in U}$, respectively. 

\subsection{MPC Policy} 
With perfect knowledge of the system parameters, MPC can be used to land the VTVL vehicle on the designated landing pad~\citep{rawlings2017mpc}. Providing an appropriate reference trajectory $r$, linearizing around an appropriately selected equilibrium point and discretizing the resulting system using zero-order-hold sampling gives the MPC control scheme 
\begin{align}\label{eq:MPC}
    \min_{x\in\mathbb{R}^{6\Tf},u\in\mathbb{R}^{3\Tf}} \:\: & \sum_{i=1}^{\Tf} |x_i-r|_Q^2 + |u_i|_R^2, \\
    \textrm{s.t.} \:\: & x_0 = \hat{x}_0, 
    \:\: x_{i+1} = Ax_i + Bu_i ,  \quad (x,u) \in (X \times Y \times \R^2 \times \Theta \times \R \times U)^{\Tf}. \nonumber
\end{align}

\noindent
Numerical simulations suggest that the performance of policy \eqref{eq:MPC} is highly sensitive on the quality of the state-space model $(A,B)$. To illustrate this point, two distinct models are considered. The first model, named model A, uses the exact parameters of model \eqref{eq:dynamics}. Model B uses inaccurately estimated parameters adjusted by increasing $l_1$ by 33\% and by decreasing $l_2$ by 25\%. Additionally, we set $m$ to be time-invariant. In model A, $m$ is set to the initial true weight $m(0)$, but is set to $\frac{1}{2}m(0)$ for model B. Fig.~\ref{fig:landing_figs} illustrates the performance of the MPC policy. The rocket successfully lands using MPC~(A), confirming the effectiveness of a well-estimated model coupled with a linearized MPC scheme \eqref{eq:MPC}. Performance significantly deteriorated with MPC~(B),  demonstrating that an inaccurate model can lead to unsatisfactory MPC performance on this system.

\subsection{DeePC-Hunt Policy} 

Two DeePC-Hunt policies are evaluated, both using the same offline training data \(w_d\) to construct the Hankel matrices in \eqref{eq:DeePC}. One policy optimizes \(\lambda\) using model A to define \(\pi^{\lambda_A}\), while the other uses model B to define \(\pi^{\lambda_B}\). Both policies employ the DeePC formulation in \eqref{eq:DeePC} with the regularization function in \eqref{eq:regulariser}, setting \(Q = \textrm{diag}(100,10,5,1,3000,30)\), \(R = \textrm{diag}(0.01,0.01,0.01)\), \(\Tf = 10\), \(\Tini = 1\), input constraints \(\mathcal{U} = U^{\Tf}\) and output constraints \(\mathcal{Y} = (X \times Y \times \mathbb{R}^2 \times \Theta \times \mathbb{R})^{\Tf}\) where \(X = [0,33.33]\), \(Y = [0,26.66]\), \(\Theta = [-0.61,0.61]\) and \(U = [0,16118.5] \times [0,322.37] \times [-0.26,0.26]\). For the DeePC-Hunt training routine described in \algorithmref{alg:deepc}, we initialize \(\lambda_0 = (50,50,1000)\), set \(N_{\text{iter}} = 100\), \(N = 20\), \(B = 1\), \((\eta^{\text{max}},\eta^{\text{min}},\beta,\alpha) = (10^2,10^{-3},1.2,0.5)\), and \(\Lambda = [10^{-5}, 10^5]^3\). Upon termination, the optimized parameters are \(\lambda_A = (49.84, 8.36, 1000.05)\) for model A and \(\lambda_B = (27.475, 2.128, 946.05)\) for model B.

\begin{figure}
    \centering
    \begin{minipage}{.2\linewidth}
        \begin{tikzpicture}
            \node[anchor=south west, inner sep=0] (image) at (0,0) {\includegraphics[trim={1.25cm} {1.25cm} {1.25cm} {1.25cm},clip,width=\linewidth]{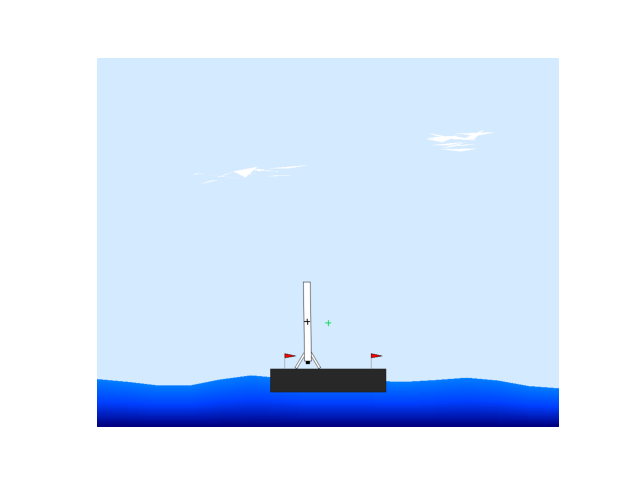}};
            \node[anchor=south] at ($(image.north) + (0, 0.1)$) {\text{MPC (A)}};
        \end{tikzpicture}
    \end{minipage}
    \hfill
    \begin{minipage}{.28\linewidth}
        \centering
        \includegraphics[width=\linewidth]{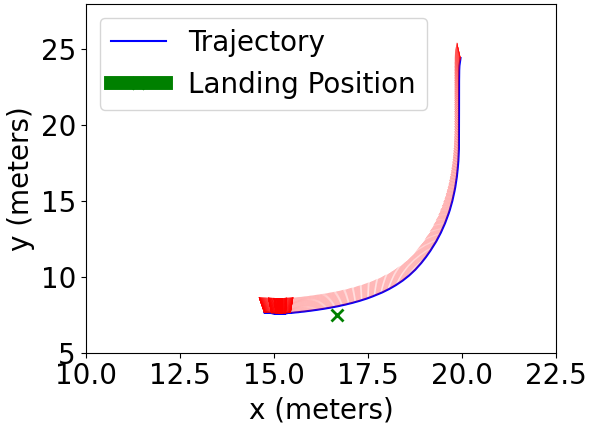}
    \end{minipage}
    \hfill
    \begin{minipage}{.2\linewidth}
        \begin{tikzpicture}
            \node[anchor=south west, inner sep=0] (image) at (0,0) {\includegraphics[trim={1.25cm} {1.25cm} {1.25cm} {1.25cm},clip,width=\linewidth]{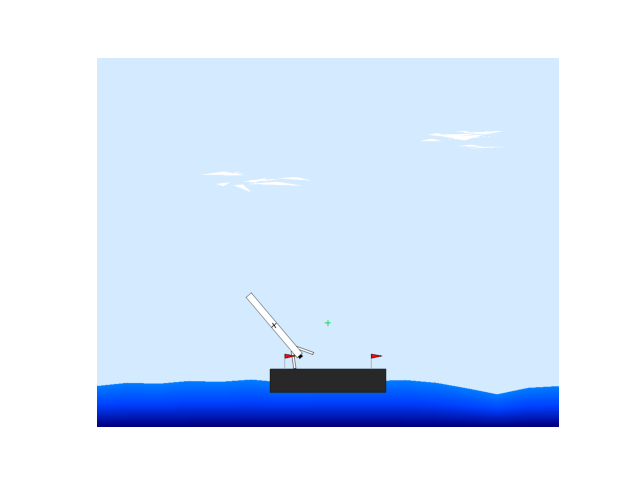}};
            \node[anchor=south] at ($(image.north) + (0, 0.1)$) {\text{MPC (B)}};
        \end{tikzpicture}
    \end{minipage}
    \hfill
    \begin{minipage}{.28\linewidth}
        \centering
        \includegraphics[width=\linewidth]{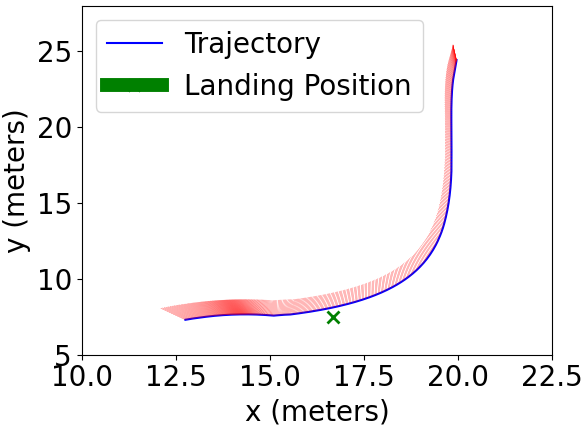}
    \end{minipage}
    \vskip\baselineskip
    \begin{minipage}{.2\linewidth}
        \begin{tikzpicture}
            \node[anchor=south west, inner sep=0] (image) at (0,0) {\includegraphics[trim={1.25cm} {1.25cm} {1.25cm} {1.25cm},clip,width=\linewidth]{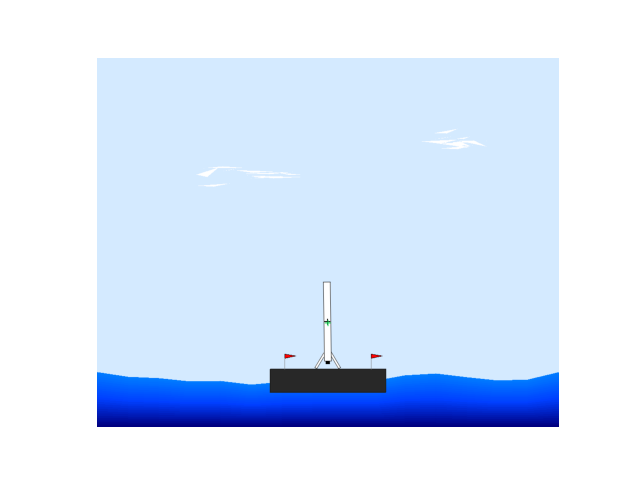}};
            \node[anchor=south] at ($(image.north) + (0, 0.1)$) {\text{DeePC-Hunt (A)}};
        \end{tikzpicture}
    \end{minipage}
    \hfill
    \begin{minipage}{.28\linewidth}
        \centering
        \includegraphics[width=\linewidth]{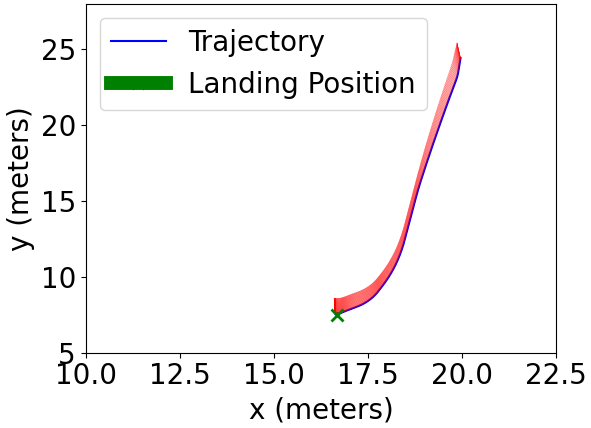}
    \end{minipage}
    \hfill
    \begin{minipage}{.2\linewidth}
        \begin{tikzpicture}
            \node[anchor=south west, inner sep=0] (image) at (0,0) {\includegraphics[trim={1.25cm} {1.25cm} {1.25cm} {1.25cm},clip,width=\linewidth]{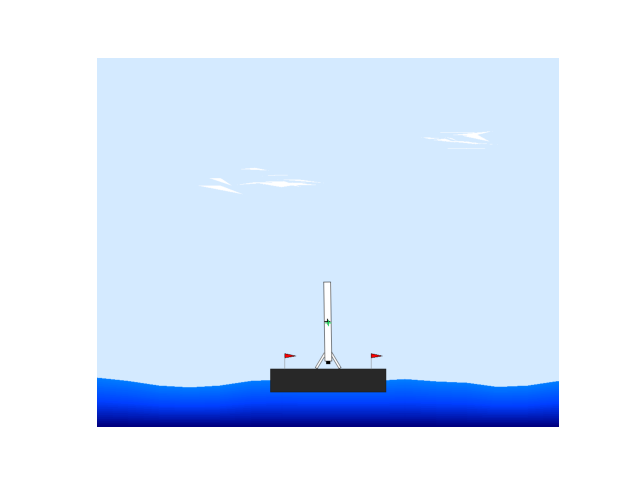}};
            \node[anchor=south] at ($(image.north) + (0, 0.1)$) {\text{DeePC-Hunt (B)}};
        \end{tikzpicture}
    \end{minipage}
    \hfill
    \begin{minipage}{.28\linewidth}
        \centering
        \includegraphics[width=\linewidth]{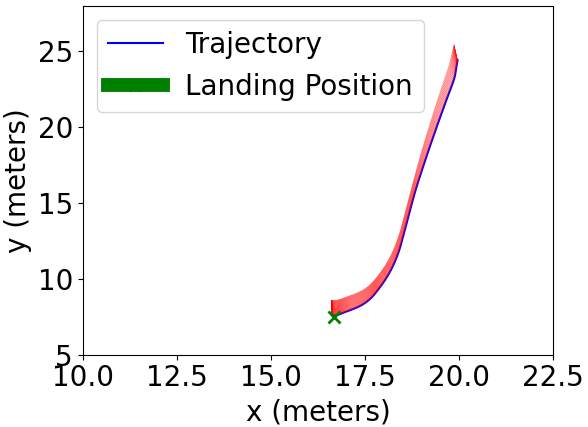}
    \end{minipage}
    \caption{Performance of the MPC (top) and DeePC-Hunt (bottom) policies using models A (left) and B (right), respectively. }
    \label{fig:landing_figs}
\end{figure}
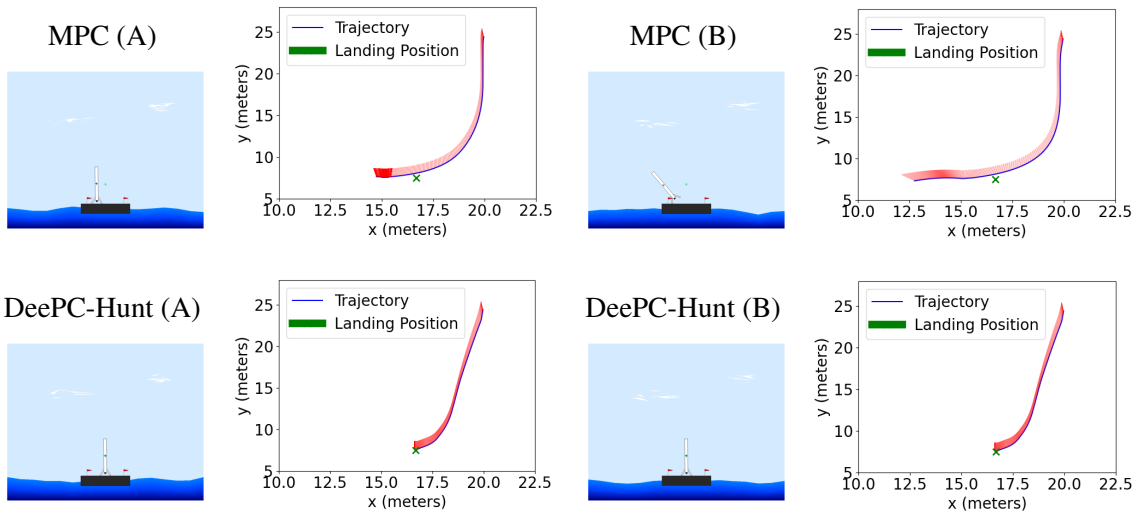

Training data is collected by applying a persistently exciting input trajectory \(u_{\text{d}} \in \mathbb{R}^{mT}\) of length \(T \in \mathbb{N}\) to the VTVL vehicle and recording the corresponding output trajectory \(y_{\text{d}} \in \mathbb{R}^{pT}\). The input trajectory is a Pseudo-Random Binary Sequence (PRBS) sequence starting from an equilibrium point \(y_{\text{eq}} \in \mathbb{R}^p\). Figure \ref{fig:landing_figs} 
illustrates that both DeePC-Hunt policies demonstrate nearly identical trajectories, even with model inaccuracies, highlighting the robustness of DeePC-Hunt. 

\subsection{Evaluation of Control Performance}

We evaluate control performance using two metrics: (i) the realized cost, measured as the total system cost incurred over \(N_\pi\) time steps taken by policy \(\pi\) to land the VTVL vehicle, and (ii) the empirical success rate across \(k = 50\) initial states uniformly sampled from a predefined range \((x, y) \in X_s \times Y_s\). The success rate represents the percentage of successful landings from these initial conditions. Results for MPC and DeePC-Hunt policies are summarized in Table \ref{tab:deepc_results}.

\begin{wraptable}{r}{0.42\textwidth}
\vspace{-10pt}
    \centering
    \scalebox{0.9}{
    \begin{tabular}{|c|c|c|}
    \hline
     & \textbf{Cost ($\times 10^3$)} & \textbf{\% Success}\\
    \hline
        \text{MPC (A)} & \text{11.228}  & 46\%\\
        $\pi^{\lambda_A}$ & 16.678 & \text{56\%}\\
    \hline
        \text{MPC (B)} & \text{8.213} & 22\% \\
        $\pi^{\lambda_B}$ & 12.165 & \text{66\%} \\
    \hline
    \end{tabular}
    }
    \vspace{-5pt}
    \caption{Results for DeePC-Hunt and MPC policies using models A and B,
    with ${X_s = [6.67,26.66]}$ and ${Y_s = [18.66,24]}$.}
\label{tab:deepc_results}
\vspace{-10pt}
\end{wraptable}

\noindent
For the sake of a fair comparison, we note that both DeePC-Hunt and MPC leverage approximate linear models, but their training pipelines differ: DeePC-Hunt relies on the system dynamics \eqref{eq:model_approximate} for optimizing $\lambda$, while MPC uses a fixed linearization of \eqref{eq:model_approximate}.
Despite this asymmetry, DeePC-Hunt yields robust closed-loop performance and 
significantly outperforms MPC in success rate in the presence of model inaccuracies. However, a trade-off emerges as MPC exhibits lower closed-loop costs upon successful landings. This discrepancy primarily arises from the impact of regularization  on  the cost function in \eqref{eq:DeePC}, prioritizing regularization $(g,\sigma_y)$ over progress towards the setpoint $\wr$. Employing model B with DeePC-Hunt yields significantly improved performance over direct MPC.

\section{Conclusion}\label{sec:conclusion}
This paper introduced DeePC-Hunt, a backpropagation-based method for automatic hyperparameter tuning of the DeePC algorithm. We demonstrated its effectiveness, showing that DeePC-Hunt outperforms model-mismatched MPC on a challenging nonlinear control task without manual tuning. Our results highlight that even with a poor approximation of the true dynamics, effective regularization parameters can be identified using backpropagation to ensure strong closed-loop performance on the true system. Given its simplicity and ease of implementation, DeePC-Hunt offers a valuable solution for automated parameter tuning, especially when exploring multiple parameters is expensive, but obtaining an approximate model is comparatively simpler.

\acks{The authors would like to thank Dylan Vogel and Gerasimos Maltezos for contributing to the development of the Gym environment, as well as Joshua N$\ddot{\text{a}}$f for writing the code to generate the trajectory plots. }

\bibliography{main}

\begin{thebibliography}{30}
\providecommand{\natexlab}[1]{#1}
\providecommand{\url}[1]{\texttt{#1}}
\expandafter\ifx\csname urlstyle\endcsname\relax
  \providecommand{\doi}[1]{doi: #1}\else
  \providecommand{\doi}{doi: \begingroup \urlstyle{rm}\Url}\fi

\bibitem[Agrawal et~al.(2018)Agrawal, Verschueren, Diamond, and Boyd]{rewriting}
Akshay Agrawal, Robin Verschueren, Steven Diamond, and Stephen Boyd.
\newblock A rewriting system for convex optimization problems.
\newblock \emph{Journal of Control and Decision}, 5\penalty0 (1):\penalty0 42--60, 2018.

\bibitem[Agrawal et~al.(2019{\natexlab{a}})Agrawal, Amos, Barratt, Boyd, Diamond, and Kolter]{layers}
Akshay Agrawal, Brandon Amos, Shane Barratt, Stephen Boyd, Steven Diamond, and J~Zico Kolter.
\newblock Differentiable convex optimization layers.
\newblock \emph{Advances in neural information processing systems}, 32, 2019{\natexlab{a}}.

\bibitem[Agrawal et~al.(2019{\natexlab{b}})Agrawal, Barratt, Boyd, Busseti, and Moursi]{agrawal2020differentiating}
Akshay Agrawal, Shane Barratt, Stephen Boyd, Enzo Busseti, and Walaa~M Moursi.
\newblock Differentiating through a cone program.
\newblock \emph{arXiv preprint arXiv:1904.09043}, 2019{\natexlab{b}}.

\bibitem[Agrawal et~al.(2020)Agrawal, Barratt, Boyd, and Stellato]{agrawal2019learning}
Akshay Agrawal, Shane Barratt, Stephen Boyd, and Bartolomeo Stellato.
\newblock Learning convex optimization control policies.
\newblock In \emph{Learning for Dynamics and Control}, pages 361--373. PMLR, 2020.

\bibitem[Amos and Kolter(2017)]{amos2021optnet}
Brandon Amos and J~Zico Kolter.
\newblock Optnet: Differentiable optimization as a layer in neural networks.
\newblock In \emph{International conference on machine learning}, pages 136--145. PMLR, 2017.

\bibitem[Amos et~al.(2018)Amos, Jimenez, Sacks, Boots, and Kolter]{amos2019differentiable}
Brandon Amos, Ivan Jimenez, Jacob Sacks, Byron Boots, and J~Zico Kolter.
\newblock Differentiable {MPC} for end-to-end planning and control.
\newblock \emph{Advances in neural information processing systems}, 31, 2018.

\bibitem[Antoulas(2005)]{antoulas2005approximation}
Athanasios~C Antoulas.
\newblock \emph{Approximation of large-scale dynamical systems}.
\newblock SIAM, 2005.

\bibitem[Borrelli et~al.(2017)Borrelli, Bemporad, and Morari]{borelliBook}
Francesco Borrelli, Alberto Bemporad, and Manfred Morari.
\newblock \emph{Predictive Control for Linear and Hybrid Systems}.
\newblock Cambridge University Press, 2017.

\bibitem[Brockman et~al.(2016)Brockman, Cheung, Pettersson, Schneider, Schulman, Tang, and Zaremba]{brockman2016openai}
G.~Brockman, V.~Cheung, L.~Pettersson, J.~Schneider, J.~Schulman, J.~Tang, and W.~Zaremba.
\newblock {OpenAI Gym}.
\newblock \emph{arXiv preprint arXiv:1904.09043}, 2016.

\bibitem[Chiuso et~al.(2023)Chiuso, Fabris, Breschi, and Formentin]{chiuso2023harnessing}
Alessandro Chiuso, Marco Fabris, Valentina Breschi, and Simone Formentin.
\newblock Harnessing uncertainty for a separation principle in direct data-driven predictive control.
\newblock \emph{arXiv preprint arXiv:2312.14788}, 2023.

\bibitem[Coulson et~al.(2019{\natexlab{a}})Coulson, Lygeros, and D{\"o}rfler]{coulson2019dataenabled}
Jeremy Coulson, John Lygeros, and Florian D{\"o}rfler.
\newblock Data-enabled predictive control: In the shallows of the deepc.
\newblock \emph{European Control Conference (ECC)}, 2019{\natexlab{a}}.

\bibitem[Coulson et~al.(2019{\natexlab{b}})Coulson, Lygeros, and D{\"o}rfler]{coulson2019regularized}
Jeremy Coulson, John Lygeros, and Florian D{\"o}rfler.
\newblock Regularized and distributionally robust data-enabled predictive control.
\newblock \emph{Conference on Decision and Control (CDC)}, 2019{\natexlab{b}}.

\bibitem[Coulson et~al.(2021)Coulson, Lygeros, and D{\"o}rfler]{distributionally_robust_deepc}
Jeremy Coulson, John Lygeros, and Florian D{\"o}rfler.
\newblock Distributionally robust chance constrained data-enabled predictive control.
\newblock \emph{IEEE Transactions on Automatic Control}, 67\penalty0 (7):\penalty0 3289--3304, 2021.

\bibitem[Diamond and Boyd(2016)]{diamond2016cvxpy}
Steven Diamond and Stephen Boyd.
\newblock {CVXPY}: {A} {P}ython-embedded modeling language for convex optimization.
\newblock \emph{Journal of Machine Learning Research}, 17\penalty0 (83):\penalty0 1--5, 2016.

\bibitem[Dörfler et~al.(2023)Dörfler, Coulson, and Markovsky]{dörfler2021bridging}
Florian Dörfler, Jeremy Coulson, and Ivan Markovsky.
\newblock Bridging direct and indirect data-driven control formulations via regularizations and relaxations.
\newblock \emph{IEEE Transactions on Automatic Control}, 68\penalty0 (2):\penalty0 883--897, 2023.

\bibitem[Favoreel and De~Moor(1999)]{SPC}
Wouter Favoreel and Bart De~Moor.
\newblock S{P}{C}: Subspace predictive control.
\newblock \emph{IFAC Proceedings Volumes}, 32, 1999.

\bibitem[Ferrante(2017)]{Ferrante2017ARC}
Reuben Ferrante.
\newblock A robust control approach for rocket landing, 2017.

\bibitem[Hu et~al.(2023)Hu, Zhang, Li, Mesbahi, Fazel, and Ba{\c{s}}ar]{policy_optimisation}
Bin Hu, Kaiqing Zhang, Na~Li, Mehran Mesbahi, Maryam Fazel, and Tamer Ba{\c{s}}ar.
\newblock Toward a theoretical foundation of policy optimization for learning control policies.
\newblock \emph{Annual Review of Control, Robotics, and Autonomous Systems}, 6\penalty0 (1):\penalty0 123--158, 2023.

\bibitem[Parikh and Boyd(2014)]{parikh2014proximal}
Neal Parikh and Stephen Boyd.
\newblock Proximal algorithms.
\newblock \emph{Foundations and trendsin Optimization}, 1\penalty0 (3):\penalty0 127--239, 2014.

\bibitem[Paszke et~al.(2019)Paszke, Gross, Massa, Lerer, Bradbury, Chanan, Killeen, Lin, Gimelshein, Antiga, Desmaison, Kopf, Yang, DeVito, Raison, Tejani, Chilamkurthy, Steiner, Fang, Bai, and Chintala]{Torch}
Adam Paszke, Sam Gross, Francisco Massa, Adam Lerer, James Bradbury, Gregory Chanan, Trevor Killeen, Zeming Lin, Natalia Gimelshein, Luca Antiga, Alban Desmaison, Andreas Kopf, Edward Yang, Zachary DeVito, Martin Raison, Alykhan Tejani, Sasank Chilamkurthy, Benoit Steiner, Lu~Fang, Junjie Bai, and Soumith Chintala.
\newblock Pytorch: An imperative style, high-performance deep learning library.
\newblock \emph{Advances in Neural Information Processing Systems 32}, 2019.

\bibitem[Rawlings et~al.(2017)Rawlings, Mayne, Diehl, et~al.]{rawlings2017mpc}
James~Blake Rawlings, David~Q Mayne, Moritz Diehl, et~al.
\newblock \emph{Model predictive control: theory, computation, and design}, volume~2.
\newblock Nob Hill Publishing Madison, WI, 2017.

\bibitem[Riedmiller and Braun(1993)]{Rprop}
Martin Riedmiller and Heinrich Braun.
\newblock A direct adaptive method for faster backpropagation learning: the rprop algorithm.
\newblock \emph{International Conference on Neural Networks}, 1993.

\bibitem[Rumelhart et~al.(1986)Rumelhart, Hinton, and Williams]{backprop}
David~E Rumelhart, Geoffrey~E Hinton, and Ronald~J Williams.
\newblock Learning representations by back-propagating errors.
\newblock \emph{nature}, 323\penalty0 (6088):\penalty0 533--536, 1986.

\bibitem[Ryu and Yin(2022)]{LSCO}
Ernest~K. Ryu and Wotao Yin.
\newblock \emph{Large-Scale Convex Optimization: Algorithms \& Analyses via Monotone Operators}.
\newblock Cambridge University Press, 2022.

\bibitem[Song et~al.(2024)Song, Kim, and Scaramuzza]{song2024learning}
Yunlong Song, Sangbae Kim, and Davide Scaramuzza.
\newblock Learning quadruped locomotion using differentiable simulation.
\newblock \emph{arXiv preprint arXiv:2403.14864}, 2024.

\bibitem[van Waarde et~al.(2020)van Waarde, De~Persis, Camlibel, and Tesi]{persis}
Henk~J. van Waarde, Claudio De~Persis, M.~Kanat Camlibel, and Pietro Tesi.
\newblock Willems’ fundamental lemma for state-space systems and its extension to multiple datasets.
\newblock \emph{IEEE Control Systems Letters}, 4\penalty0 (3):\penalty0 602--607, 2020.

\bibitem[Waarde et~al.(2020)Waarde, Eising, Trentelman, and Camlibel]{vanwaarde2020data_informativity}
Henk Waarde, Jaap Eising, Harry Trentelman, and Kanat Camlibel.
\newblock Data informativity: A new perspective on data-driven analysis and control.
\newblock \emph{IEEE Transactions on Automatic Control}, pages 1--1, 01 2020.

\bibitem[Wiedemann et~al.(2023)Wiedemann, W{\"u}est, Loquercio, M{\"u}ller, Floreano, and Scaramuzza]{wiedemann2023training}
Nina Wiedemann, Valentin W{\"u}est, Antonio Loquercio, Matthias M{\"u}ller, Dario Floreano, and Davide Scaramuzza.
\newblock Training efficient controllers via analytic policy gradient.
\newblock \emph{International Conference on Robotics and Automation (ICRA)}, 2023.

\bibitem[Xue and Matni(2021)]{xue2021datadriven_SLS}
Anton Xue and Nikolai Matni.
\newblock Data-driven system level synthesis.
\newblock \emph{Learning for Dynamics and Control}, 2021.

\bibitem[Zuliani et~al.(2023)Zuliani, Balta, and Lygeros]{zuliani2024bpmpc}
Riccardo Zuliani, Efe~C Balta, and John Lygeros.
\newblock Bp-mpc: Optimizing closed-loop performance of mpc using backpropagation.
\newblock \emph{arXiv preprint arXiv:2312.15521}, 2023.

\end{thebibliography}

\end{document}